\documentclass[11pt]{amsart}
\input amssym.def
\usepackage{amsbsy}
\newtheorem{theorem}{Theorem}
\newtheorem{definition}{Definition}

\newtheorem*{acknowledgement*}{Acknowledgement}

\begin{document}

\title[On determination of Zero-sum $\ell$-generalized Schur Numbers]{On
determination of Zero-sum $\ell$-generalized Schur Numbers for some linear
equations}

\author{Bidisha Roy and Subha Sarkar}
   \address[Bidisha Roy and Subha Sarkar]{Harish-Chandra Research Institute, HBNI\\
Chhatnag Road, Jhunsi\\
Allahabad 211019, India}
\email[Subha Sarkar]{subhasarkar@hri.res.in}
\email[Bidisha Roy]{bidisharoy@hri.res.in}

\begin{abstract}
Let  $r$, $m$ and $k\geq 2$ be positive integers such that $r\mid k$ and let  $v \in \left[ 0,\lfloor \frac{k-1}{2r} \rfloor \right]$ be any integer. For any integer $\ell \in [1, k]$ and $\epsilon \in \{0,1\}$, we let $\mathcal{E}_{v}^{(\ell, \epsilon)}$ be the linear homogeneous equation defined by $\mathcal{E}_{v}^{(\ell, \epsilon)}: x_1 + \cdots + x_{k-(rv+\epsilon)} =x_{k-(rv+\epsilon-1)} +\cdots+ \ell x_{k}$. We denote the number $S_{\mathfrak{z},m}^{(\ell, \epsilon)}(k;r;v)$, which is defined to be the least positive integer $t$ such that for any $m$-coloring $\chi: [1, t] \to
\{0, 1,\ldots,m-1\}$, there exists a solution $(\hat{x}_1, \hat{x}_2, \ldots, \hat{x}_k)$ to the equation $\mathcal{E}_{v}^{(\ell,\epsilon)}$ that satisfies the $r$-zero-sum condition, namely, $\displaystyle\sum_{i=1}^k\chi(\hat{x}_i) \equiv 0\pmod{r}$. In this article, we completely determine the constant  $S_{\mathfrak{z}, 2}^{(k,1)}(k;r;0)$, $S_{\mathfrak{z}, m}^{(k-1,1)}(k;r;0)$, $S_{\mathfrak{z}, 2}^{(1,1)}(k;2;1)$ and $S_{\mathfrak{z}, r}^{(1,0)}(k;r;v)$. Also, we prove upper bound for the constants $S_{\mathfrak{z},2}^{(2,1)}(k;2;0)$ and $S_{\mathfrak{z},2}^{(1,1)}(k;2;v)$.
\end{abstract}
\maketitle

\section{Introduction}
For a given positive integer $m$, the classical {\it Schur number}, denoted by $S(m)$ is defined as follows. For a given positive integer $m$, the number $S(m)$ is defined to be the least positive integer $t$ such that every $m$-coloring of the interval $[1,t]$ admits a monochromatic solution to the linear equation $ x_1 +x_2 =x_3$.

\medskip

In 1933, Rado \cite{rado} introduced the {\it generalized Schur number} and defined it as follows. For given positive integers $k\geq 2$ and $m\geq 1$, the {\it generalized Schur number} $S(k;m)$ is the least positive integer $t$ such that every $m$-coloring of the interval $[1,t]$ admits a monochromatic solution to the equation $ x_1 +\cdots +x_{k-1} = x_k$. He also proved that $S(k;m)$ is finite. Moreover, he characterized any linear homogeneous equation with non-zero coefficients that admits a monochromatic solutions as follows.

\begin{theorem}\cite{rado}\label{rado}
Let $\mathcal{E}: a_1x_1+\cdots + a_kx_k = 0$ be any linear homogeneous equation with integer coefficients and at least one of the coefficients is non-zero. Then $\mathcal{E}$ admits monochromatic solution to any $m$-coloring of natural numbers  if and only if there exist $1\leq i_1 < i_2 < \ldots < i_d \leq k$ such that $a_{i_1}+a_{i_2}+\cdots+a_{i_d} = 0$. 
\end{theorem} 

Before we proceed further, we need the following definition.

\begin{definition}
For a given positive integer $r\geq 2$, we call a sequence $(a_1,a_2,\ldots,a_n)$ of integers to be an {\it $r$-zero-sum} sequence if
 $\sum\limits_{i=1}^{n} a_i \equiv 0 \pmod r$.
\end{definition}

The zero-sum theory has received a lots of attention and has vast literature, see for
instance \cite{bd}, \cite{caro}, \cite{egz} and \cite{gg}.

\medskip

In \cite{aaron}, Robertson replaced the ``monochromatic property"  in the definition of the {\it generalized Schur number} by ``zero-sum property" and introduced  {\it zero-sum generalized Schur number}. A similar kind of topic has been addressed in \cite{brown}.

\medskip

Let $r$ and $k$ be given positive integers such that $r \mid k$. For given integers $\ell\in [1, k]$, $v \in \left[ 0,\lfloor \frac{k-1}{2r} \rfloor \right]$ and $ \epsilon \in \{ 0,1 \}$ we set the linear homogeneous equation as follows.
\begin{equation}\label{homo-equation}
\mathcal{E}_{v}^{(\ell, \epsilon)} : x_1 + \cdots + x_{k-(rv+\epsilon)} =x_{k-(rv+\epsilon-1)} +\cdots+ \ell x_{k}
\end{equation} .

\begin{definition} (Zero sum $\ell$-generalized Schur number)\label{def}
Let $r$, $m$ and $k\geq 2$ be given positive integers with $r|k$ and let $v \in \left[ 0,\lfloor \frac{k-1}{2r} \rfloor \right]$ be any integer and $\epsilon\in\{0,1\}$. For any integer $\ell \in [1, k]$,  the {\it zero-sum $\ell$-generalized Schur number} $S_{\mathfrak{z},m}^{(\ell, \epsilon)}(k;r;v)$ is defined to be the least positive integer $t$ such that for every $m$-coloring $\chi :[1,t] \rightarrow \{0,1,\ldots,m-1\}$, there is a solution $(x_1,x_2,\ldots,x_k)$ to the equation $\mathcal{E}_{v}^{(\ell, \epsilon)}$ satisfying $r$-zero-sum condition, namely,  $\sum\limits_{i=1}^ {k} \chi ( x_i) \equiv 0\pmod r $.
\end{definition}

For any integers $ \ell \in [1, k]$, $v \in \left[ 1,\lfloor \frac{k-1}{2r} \rfloor \right]$ and $ \epsilon \in \{ 0,1 \}$ or $ \ell \in [1, k-1], v=0$ with $\epsilon=1$, the linear equation
$\mathcal{E}_{v}^{(\ell, \epsilon)}$ defined in \eqref{homo-equation} satisfies the condition of Theorem \ref{rado} and hence $S_{\mathfrak{z}, m}^{(\ell, \epsilon)}(k;r;v)$ are finite numbers for any $k, r$ and $m$. 

\medskip

When $\ell =\epsilon=1$, $v =0$ and $m = r$,  the constant $S_{\mathfrak{z},r} ^{(1,1)}(k;r;0)$ is denoted by $S_{\mathfrak{z}}(k;r)$ in the literature (see \cite{aaron}) and is called the {\it zero-sum generalized Schur number}. In  \cite{aaron}, Robertson proved a lower bound for $S_{\mathfrak{z}}(k,r)$ when $ r=3$ and $4$.

\medskip

Recently in \cite{erik}, E. Metz showed the exact values of this constant for $r=3 ,4$. Moreover, he proved the following two results.

\begin{theorem}\cite{erik}
Let $r$ and $k$ be positive integers with $ r| k $ and $k \geq 2r$. Then,
$$
S_{\mathfrak{z}}(k;r)=S_{\mathfrak{z},r}^{(1,1)}(k;r;0) \leq
\begin{cases}
 kr-r  & ;r \mbox{ is an odd prime }\\
 4k-5  &   ; r=4\\
 kr- \sum_{i=1}^t (p_i-1)-1  & ; r \geq 6 \mbox{ and }r = p_1 \ldots p_t \\& \mbox{ be the prime decomposition of }r\\ &\mbox{ and } p_i \mbox{'s are not necessaryly distinct }.
\end{cases}
$$         
\end{theorem}

\begin{theorem}\cite{erik}
Let $r$ and $k$ be positive integers with $ r| k $. Then,
$$
S_{\mathfrak{z}}(k;r)=S_{\mathfrak{z},r}^{(1,1)}(k;r;0) \geq
\begin{cases}
 kr-r  & ;r \mbox{ is an odd integer}\\
 kr-r-1  & ;r \mbox{ is an even integer}.
\end{cases}
$$         
\end{theorem}

Also, when  $\ell =\epsilon =1$, $ v=0$ and $m=2$ in the Definition \ref{def}, the constant $S_{\mathfrak{z},2}^{(1,1)}(k;r;0)$ is denoted by $S_{\mathfrak{z},2}(k;r)$ in the literature. In \cite{aaron}, Robertson calculated the value of this constant for some cases and was completely determined in \cite{abs}. More precisely,

\begin{theorem}\cite{aaron,abs}
Let $r$ and $k$ be positive integers with $ r| k $ and $k > r$. Then,
$$
S_{\mathfrak{z},2}(k;r)=S_{\mathfrak{z},2}^{(1,1)}(k;r;0) = rk-2r+1.
$$         
\end{theorem}

In this article, we shall consider the case $\ell >1$. First note that, when $\ell = k-1, v = 0$ and $ \epsilon =1$, the condition of Theorem \ref{rado} is satisfied by $\mathcal{E}_{0}^{(k-1,1)}$. Also $(1,1,\ldots,1)$ satisfies the equation $\mathcal{E}_{0}^{(k-1,1)}$ together with the $r$-zero-sum condition. Thus we see that  $S_{\mathfrak{z},m}^{(k-1,1)}(k;r;0) = 1$ for each $m \geq 2$.
\medskip

However, the equation $\mathcal{E}_{0}^{(k,1)}: x_1+ \cdots +x_{k-1}=kx_{k}$ does not satisfy the condition of Theorem \ref{rado}. Hence, by Theorem \ref{rado}, we can not conclude the finiteness of the constant $S_{\mathfrak{z},2}^{(k,1)}(k;r;0)$. We indeed prove the finiteness of this constant by calculating the exact value of it as follows.

\begin{theorem}\label{k}
Let $r$ and $k$ be positive integers with $ r| k$ with $k\geq 2$. Then,
$$
S_{\mathfrak{z},2}^{(k,1)}(k;r;0) =
\begin{cases}
 3  & ;r=2 \text{ and } k \geq 4 \\
 4  & ;r\geq 3 \text{ and } k=r \text{ or } k=2r \\
 3  & ;r\geq 3 \text{ and } k\geq 3r.
\end{cases}
$$                           
\end{theorem}
When $\ell = m = r = 2$, $\epsilon = 1$ and $v=0$, we have the following upper bound.

\begin{theorem}\label{2}
Let $k\geq 4$ be an even positive integer. Then,
$$
S_{\mathfrak{z},2}^{(2,1)}(k;2;0) \leq \left\lceil \frac{k}{4}\right\rceil+\frac{k}{2}-1.
$$            
\end{theorem}

Now we move for the case when $v$ is not zero in the equation $\mathcal{E}_{v}^{(\ell, \epsilon)}$.

\begin{theorem}\label{lower bound}
Let $k$ be an even positive integer and $v \in [1, \left\lfloor \frac{k-1}{4} \right\rfloor]$. Then $$S_{\mathfrak{z},2}^{(1,1)}(k;2;v) \leq \left(\frac{k}{2}-2v\right).$$
\end{theorem}

When $v=1$ in Theorem \ref{lower bound}, we prove the exact value in the following theorem. 
\begin{theorem}\label{general}
%

Let $k$ be an even positive integer with $k \geq 6$. Then, 
$$S_{\mathfrak{z},2}^{(1,1)}(k;2;1) = \left(\frac{k}{2} -u-2\right), \mbox{ where } u =
\begin{cases}
 t  & ;\text{ if } k=10t+s \text{ and } s \in \{6,8\} \\
 t-1  & ;\text{ if } k=10t+s \text{ and } s \in \{0,2,4\}
\end{cases}.$$
\end{theorem}

In the following theorem, we compute the exact value when $\epsilon = 0$, $l= 1$ and $m = r$.

\begin{theorem}\label{more}
Let  $r$ and $k$ are positive integers with $ r| k $ and $v \in [1, \left\lfloor \frac{k-1}{2r} \right\rfloor]$ also an integer. Then, 
$$
S^{(1,0)}_{\mathfrak{z},r}(k;r;v)= \frac{k}{r} - \left\lfloor \frac{(v-1)k}{vr} \right\rfloor -1.
$$
\end{theorem}

\medskip

\section{proof of the theorem \ref{k}}

\noindent Here we consider the equation $\mathcal{E}_{0}^{(k,1)}: x_1+ \cdots +x_{k-1}=kx_{k}.$

\medskip

\noindent{\sf Case I:} ($r=2$ and $k \geq 4$) 

\medskip

Let  $\chi: [1,3] \rightarrow \{0,1\}$ be any $2$-coloring.  We may assume that
$\chi(1) = 0$, since $\chi$ admits a $2$-zero sum solution if and only if $\hat{\chi}$
defined by $\hat{\chi}(i) = 1- \chi(i)$ does. Considering the solution
$(\underbrace{1, \ldots, 1}_{\textrm{(k-2)-times}},2,1)$ and using the color of $1$, we can conclude $\chi(2)=1$ (for otherwise, we are done).  Since $2 \chi(3) + (k-3) \chi(2)+ \chi(2) \equiv 0 \pmod 2,$ we  see that  $( 3,3, \underbrace{2, \ldots,2}_{\textrm{(k-3)-times}},2)$ is a $2$-zero sum solution to the equation
$\mathcal{E}_{0}^{(k,1)}$. This proves that, in this case, we get  $S_{\mathfrak{z},2}^{(k,1)}(k;2;0) \leq 3$.

\medskip

For proving the lower bound, we consider the coloring $\chi: [1,2] \rightarrow \{0,1\}$ with $\chi(1)=0$ and $\chi(2)=1$.  If $(x_1, \ldots, x_k)$ is a solution to $\mathcal{E}_{0}^{(k,1)}$, then  $x_k \ne 2$ because $x_1+\cdots+x_{k-1} \leq 2(k-1)$ where as $kx_k = 2k$. If $x_k =1$, then $kx_k = k$ and hence  we can get only one solution to the equation
$\mathcal{E}_{0}^{(k,1)}$ which is $(\underbrace{1, \ldots,
1}_{\textrm{(k-2)-times}},2,1)$ but it does not satisfy $2$-zero-sum condition.  Hence, we can conclude $S_{\mathfrak{z},2}^{(k,1)}(k;2;0)= 3$.

\bigskip

\noindent{\sf{Case II:}} $(r \geq 3$  and $k=r$)  or $(r \geq 3$ and $k=2r$) 

\medskip

Let $\chi: [1,4] \rightarrow \{0,1\}$ be any $2$-coloring. We may assume that $\chi(1) = 0$. Looking at the solution $(\underbrace{1, \ldots, 1}_{\textrm{(k-2)-times}},2,1)$ and using the color of $1$, we can conclude $\chi(2)=1$ (for otherwise, we are done). Now, considering the solution
$(\underbrace{2, \ldots, 2}_{\textrm{(k-3)-times}},3,3,2)$ and using the color of $2$, we can conclude $\chi(3)=0$ (for otherwise, we are done). Again, consider the solution $(\underbrace{2, \ldots, 2}_{\textrm{(k-2)-times}},4,2)$. Then, if $\chi(4) = 1$, then we are done.  If $\chi(4)=0$, then by observing that $(\underbrace{3, \ldots, 3}_{\textrm{(k-4)-times}},4,4,4,3)$ is an $r$-zero-sum solution to the equation $\mathcal{E}_{0}^{(k,1)}$. Hence, since $r\geq 3$ and $k$ is any integer such that $r|k$,  we get $S_{\mathfrak{z},2}^{(k,1)}(k;r;0) \leq 4$.

\smallskip

Now it remains to prove the lower bound. We consider $\chi: [1, 3] \rightarrow \{0, 1\}$ such that $\chi(1) = 0 = \chi(3)$ and $\chi(2) = 1$.  

\bigskip

\noindent{\sf{Subcase 1.}} ($r\geq 3$ and $k=r$)

\bigskip

Since $k=r$, first one can observe that an $r$-zero sum solution to the equation $\mathcal{E}_{0}^{(k,1)}$  is also  a monochromatic solution  and vice versa. 

Second, note that by taking  $x_i =2$ for all $i = 1, 2,\ldots, k$, we cannot get any solution to the equation $\mathcal{E}_{0}^{(k,1)}$. Also, since $kx_k \geq r$, by taking only $x_i = 3$ or $x_i = 1$, one
can not get any solution to the equation $\mathcal{E}_{0}^{(k,1)}$. Thus,  under this coloring,  it is impossible to get a monochromatic solution to the equation $\mathcal{E}_{0}^{(k,1)}$. Hence, we can conclude
$S_{\mathfrak{z},2}^{(k,1)}(k;r;0) \geq 4$. 

\medskip

\noindent{\sf{Subcase 2.}} ($r\geq 3$ and $k=2r$) 

\medskip

By taking $x_k = 1$, the only possible solution
is $(\underbrace{1, \ldots, 1}_{\textrm{(k-2)-times}},2,1)$, which is not an $r$-zero sum solution. Again, by taking $x_k =2$, the only possible solution to the equation $\mathcal{E}_{0}^{(k,1)}$ is $(3,3, \underbrace{2, \ldots,2}_{\textrm{(k-3)-times}},2)$, 
which is not an $r$-zero sum solution under the coloring $\chi$. 
Finally, if we take $x_k =3$, then, $kx_k=3k$ and $x_1+\cdots+x_{k-1} \leq 3k-3$. Hence, there is no solution with $x_k =3$ also. Thus we get,  $S_{\mathfrak{z},2}^{(k,1)}(k;r;0) \geq 4$.

\medskip

\noindent{\sf{Case III:}} ($r\geq 3$ and $k \geq 3r$)

\medskip

Suppose for a contradiction that there exists a
$2$-coloring $\chi: [1,3] \rightarrow \{0,1\}$ for which $\mathcal{E}_{0}^{(k,1)}$ does not have any $r$-zero-sum solution for some $r \geq 3$. We may assume that $\chi(1) = 0$. Looking at the solution $(\underbrace{1, \ldots, 1}_{\textrm{(k-2)-times}},2,1)$ and
using the color of $1$, we can conclude $\chi(2)=1$. Now, considering the solution
$(\underbrace{2, \ldots, 2}_{\textrm{(k-3)-times}},3,3,2)$ and using the color of
$2$, we can conclude $\chi(3)=0$.

\medskip

If $k$ is an odd multiple of $r$, then we first observe that $\frac{k-r}{2}$ is a positive integer for any parity of $r$. Hence, in this case, we see that 
 $
 (\underbrace{2, \ldots, 2}_{\textrm{$(r-1)$-times}},\underbrace{1, \ldots,
1}_{\textrm{$\left(\frac{k-r}{2} -1\right)$-times}},\underbrace{3, \ldots,
3}_{\textrm{$\left(\frac{k-r}{2} +1\right)$-times}},2)
$ 
is an $r$-zero sum solution to the equation $\mathcal{E}_{0}^{(k,1)}$.

\medskip

If $k$ is an even multiple of $r$, then $\frac{k-2r}{2}$ is a positive
integer as  $k \geq 3r$. Thus, in this case, we see that  $(\underbrace{2, \ldots, 2}_{\textrm{$(2r-1)$-times}},\underbrace{1, \ldots,
1}_{\textrm{($\frac{k-2r}{2} -1$)-times}},\underbrace{3, \ldots,
3}_{\textrm{($\frac{k-2r}{2} +1$)-times}},2)$ is an $r$-zero sum solution to the
equation $\mathcal{E}_{0}^{(k,1)}$, which is a contradiction. Hence, we have $S_{\mathfrak{z},2}^{(k,1)}(k;r;0) \leq 3$.

\smallskip

For the lower bound, we consider the coloring $\chi: [1,2] \rightarrow \{0,1\}$
such that $\chi(1)=0$ and $\chi(2)=1$. Note that $x_k \ne 2$ as $kx_k  = 2k$ and $x_1+\cdots+x_{k-1} \leq 2k-2$. Therefore, $x_k = 1$ and hence $kx_k = k$. Thus,  there  is only one solution, namely, 
$(\underbrace{1, \ldots, 1}_{\textrm{(k-2)-times}},2,1)$, to the equation
$\mathcal{E}_{0}^{(k,1)}$, which does not satisfy $r$-zero-sum condition. This proves the lower bound and the theorem.
\hfill\(\Box\)

\section{proof of the theorem \ref{2}}
Here we consider the equation $\mathcal{E}_{0}^{(2,1)}: x_1+ \cdots +x_{k-1}=2x_{k}.$

Let us denote $S = \displaystyle\left\lceil \frac{k}{4}\right\rceil+\frac{k}{2}-1$. 
Let  $\chi: [1,\left\lceil \frac{k}{4}\right\rceil+\frac{k}{2}-1] \rightarrow \{0,1\}$ be a $2$-coloring for which
$\mathcal{E}_{0}^{(2,1)}$ does not have any $2$-zero-sum solution. We may assume that
$\chi(1) = 0$, since $\chi$ admits a $2$-zero sum solution if and only if $\hat{\chi}$
defined by $\hat{\chi}(i) = 1- \chi(i)$ does. 

\smallskip

Note that,  $\chi(\frac{k}{2})$ is either $0$ or $1$.  In the following table,  we determine the color of $2, S, \frac{k}{2}+1$ and
$3$ using some solutions of $\mathcal{E}_{0}^{(2,1)}$.
$$
\begin{array}{|l|l|l|}
\hline
\mbox{ Solution to the equation } \mathcal{E}_{0}^{(2,1)} &
\chi(\frac{k}{2})=0  & \chi(\frac{k}{2})=1\\ 
\hline
\hline
(\underbrace{1, \ldots, 1}_{\textrm{(k-2)-times}},2,\frac{k}{2}) & \chi(2)=1 &
\chi(2)=0\\\hline
(\underbrace{1, \ldots, 1}_{\textrm{(k-2)-times}},\frac{k}{2},S) & \chi(S)=1 &
\chi(S)=0\\\hline
(\underbrace{1, \ldots, 1}_{\textrm{(k-4)-times}},2,2,2,\frac{k}{2}+1) &
\chi(\frac{k}{2}+1)=0 & \chi(\frac{k}{2}+1) =1\\\hline
(\underbrace{1, \ldots, 1}_{\textrm{(k-3)-times}},2,3,\frac{k}{2}+1) & \chi(3)= 0 &
\chi(3)=0\\
\hline
\end{array}
$$

\bigskip
\noindent{\sf Case I: } $(4 | k)$ 

\bigskip

When $k=4$, we see that $(1,1,2,2)$ is a 2-zero sum
solution to the equation $x_1+x_2+x_3=2x_4$. Thus $S_{\mathfrak{z},2}^{(2,1)}(4;2;0)
\leq 2$ which is equal to  $S = \left\lceil
\frac{4}{4}\right\rceil+\frac{4}{2}-1$.

\smallskip

Let $k=4t$ for some integer $ t \geq 2$.  In this case,  $S=\left\lceil \frac{4t}{4} \right\rceil+\frac{4t}{2}-1 = 3t-1$.
Thus, using above table and considering some more solution of the equation
$\mathcal{E}_{0}^{(2,1)}$, we have the following table.
$$
\begin{array}{|l|l|l|}
\hline
\mbox{ Solution to the equation }  \mathcal{E}_{0}^{(2,1)} &
\chi(\frac{k}{2})=0 & \chi(\frac{k}{2})=1 \\
\hline
\hline
(\underbrace{1, \ldots, 1}_{\textrm{(2t+2)-times}},\underbrace{2, \ldots,
2}_{\textrm{(2t-3)-times}}, S-1) & \chi(S-1)=0 & \chi(S-1)=1\\\hline
(\underbrace{1, \ldots, 1}_{\textrm{(4t-2)-times}},2t-2, S-1) & \chi(2t-2) =1 &
\chi(2t-2)=0\\
\hline
\end{array}
$$
Thus, we get a solution $(\underbrace{1, \ldots,
1}_{\textrm{(4t-3)-times}},3,2t-2, S)$ to the equation $\mathcal{E}_{0}^{(2,1)}$ which
satisfies $2$-zero-sum condition in both cases $\chi(\frac{k}{2}) =0$ or $1$.

\bigskip

\noindent{\sf Case II: } ($4$ does not divide $k$)

\bigskip

 In this case,  we can write
$k=4t+2$ for some integer  $t \geq 1$. Also we observe that  $S=\left\lceil \frac{4t+2}{4}
\right\rceil+\frac{4t+2}{2}-1 = 3t+1$.
Thus, using the first table and considering some more solutions to the equation
$\mathcal{E}_{0}^{(2,1)}$, we get the following table;
$$
\begin{array}{|l|l|l|}
\hline
\mbox{ Solution to the equation } \mathcal{E}_{0}^{(2,1)} & 
\chi(\frac{k}{2})=0 & \chi(\frac{k}{2})=1\\
\hline
\hline
(\underbrace{1, \ldots, 1}_{\textrm{(2t+2)-times}},\underbrace{2, \ldots,
2}_{\textrm{(2t-1)-times}}, S-1) & \chi(S-1)=0 & \chi(S-1)=1\\\hline
(\underbrace{1, \ldots, 1}_{\textrm{4t-times}},2t, S-1) & \chi(2t) =1 & \chi(2t)=0\\
\hline
\end{array}
$$
Thus, we get a solution $(\underbrace{1, \ldots, 1}_{\textrm{(4t-1)-times}},3,2t,
S)$ to the equation $\mathcal{E}_{0}^{(2,1)}$ which satisfies $2$-zero-sum condition in both
cases $\chi(\frac{k}{2}) =0$ or $1$. This proves the theorem.
\hfill\(\Box\)

\medskip

\section{proof of the theorem \ref{lower bound}}

Here we consider the equation $\mathcal{E}_{v}^{(1,1)}: x_1+ \cdots +x_{k-(2v+1)}=x_{k-2v}+\cdots+x_{k}$.

\medskip

Let us denote $t : = \left( \frac{k}{2}-2v\right)$.
Let  $ \chi  :[1, t] \rightarrow \{ 0,1 \}$ be any $2$-coloring. Since $v \leq \frac{k-1}{4}$ and $k$ is even, we get $k \geq 4v+2$.

\smallskip

Now, observe that 
$$\left(\underbrace{1,\ldots,1}_{(k-4v)-\mbox{ times}} ,\underbrace{\left(\frac{k}{2}-2v\right),\ldots, \left(\frac{k}{2}-2v\right)}_{4v-\mbox{times}}\right)
$$
 is a solution to the equation 
\begin{equation}\label{same1}
x_1 + \cdots +  x_{k-(2v+1)} = x_{k-2v}+ \cdots +x_k,
\end{equation}
because $x_{k-2v}+ \cdots +x_k= (2v+1)(\frac{k}{2}-2v) = vk +\frac{k}{2}-4v^2 -2v$ and
\begin{align*}
x_1 + \ldots +x_{k-(2v+1)} & = \underbrace{1+ \ldots +1}_{(k-4v)-\mbox{times}} + \underbrace{\left(\frac{k}{2}-2v\right) + \ldots + \left(\frac{k}{2}-2v\right)}_{(2v-1)-\mbox{times}} \\
                           & = k-4v + (2v-1) (\frac{k}{2}-2v)\\
                           & = vk +\frac{k}{2}-4v^2 -2v.
                           \end{align*} 
Since $2|k$, we get
\begin{equation*}
\sum _{i=1}^ k \chi(x_{i})= (k-4v) \chi (1) + 4v \chi \left(\frac{k}{2}-2v\right) \equiv 0\pmod 2.
\end{equation*} 
Thus, we get $S_{\mathfrak{z},2}^{(1,1)}(k;2;v) \leq (\frac{k}{2}-2v).$
\hfill\(\Box\)

\medskip

\section{proof of the theorem \ref{general}}

Here we consider the equation
$\mathcal{E}_{1}^{(1,1)}:x_1 + \cdots +  x_{k-3} = x_{k-2}+ x_{k-1} +x_k$ and denote it by $\mathcal{E}(k)$ for simplicity.

\bigskip

\noindent{\sf Case I:}  ($u=0$)

\bigskip

In this case, $k \in [6, 14]$ and we show that $S_{\mathfrak{z},2}^{(1,1)}(k;2;1) = (\frac{k}{2} -2)$.

\smallskip

In this case, we first prove the lower bound.
 
 \bigskip
 
\noindent {\sf Subcase 1.}  ($k=6$).

\smallskip

In this case, clearly, $( 1,1,1,1,1,1)$ is the only $2$-zero sum solution and we are done.

\smallskip
 
\noindent{\sf Subcase 2.}  ($k=8$)

\smallskip

The lower bound follows $S_{\mathfrak{z},2}^{(1,1)}(8;2;1) \geq 2$ as $( \underbrace{1, \ldots, 1}_{\text{8-times}})$ is not a solution of $ \mathcal{E}(8)$.

\smallskip
 
\noindent{\sf Subcase 3.}  ($k=10$) 

\smallskip

$S_{\mathfrak{z},2}^{(1,1)}(10;2;1) \geq 3$ follows because $ \mathcal{E}(10)$ does not have any solution in $[1,2]$.

\smallskip
 
\noindent{\sf Subcase 4.} ($k=12$)

\smallskip

We observe that $ \mathcal{E}(12)$ does not have any solution in $[1,2]$ and the only solution of $\mathcal{E}(12)$ in $[1,3]$ is $(1, \ldots , 1 , 3 , 3, 3)$. Thus, considering the $2$-coloring $ \chi(1)= \chi(2) =0, \chi(3) =1$ of $[1,3]$, we see that $\mathcal{E}(12)$ does not have a $2$-zero sum solution in $[1,3]$. Hence, we get $S_{\mathfrak{z},2}^{(1,1)}(12;2;1) \geq 4$.

\smallskip
 
\noindent{\sf Subcase 5.} ($k=14$)

\smallskip

The only solutions of $\mathcal{E}(14)$ in $[1,4]$ are $(1, \ldots , 1, 4, 4, 3)$ and $ (1, \ldots , 1, 2, 4, 4, 4)$. Thus, if we consider the $2$-coloring $\chi(1) = \chi(2) =0, \chi(3) = \chi(4) =1 $ of $[1,4]$, then $\mathcal{E}(14)$ does not have a $2$-zero sum solution in $[1,4]$. Hence, $S_{\mathfrak{z},2}^{(1,1)}(14;2;1) \geq 5$.

\smallskip

By putting $v = 1$ in  Theorem \ref{lower bound},  the upper bound follows.

\bigskip



\noindent{\sf Case II:} ($u$ is odd)

\bigskip

For proving the upper bound, on the contrary, let $\chi : [1,\frac{k}{2}-u-2] \rightarrow \{0,1\}$ be a $2$-coloring such that $\chi$ doesn't admit a $2$-zero sum solution to $\mathcal{E}(k)$. We may assume that $\chi(1) = 0$, since $\chi$ admits a $2$-zero sum solution if and only if the induced coloring $\hat{\chi}$ defined by $\hat{\chi}(i) = 1- \chi(i)$ does.

\bigskip

\noindent{\sf Subcase 1.}  $(4 | k)$

\bigskip

Since $u$ is odd,  $u-1$ and $3u-1$ both are even. Using $\chi(1)=0$ and considering the solution $(1, \ldots, 1,\frac{k}{4}+\frac{u-1}{2}, \frac{k}{4}+\frac{u-1}{2}, \frac{k}{2}-u-2)$,  we get $\chi(\frac{k}{2}-u-2)=1$ (otherwise, we are done). Since 
$$ 
\frac{k}{2}-u-2 \geq \frac{k}{4}+\frac{3u-1}{2} \Leftrightarrow k \geq 10u+6,
$$ we can consider the solution $(1, \ldots, 1, \frac{k}{4}+\frac{3u-1}{2}, \frac{k}{4}+\frac{3u-1}{2}, \frac{k}{2}-3u-2)$ which produces $\chi(\frac{k}{2}-3u-2)=1$. But the solution $(\underbrace{1,\ldots,1}_{\text{(k-4)-times}},\frac{k}{2}-3u-2,\frac{k}{2}-u-2,\frac{k}{2}-u-2,\frac{k}{2}-u-2)$ satisfies $2$-zero sum condition and proves this subcase. 

\bigskip

\noindent{\sf Subcase 2.} ($4 \not|\ k$)

\bigskip

 In this subcase, we write $k =4t+2$  for some integer $t$, as $k$ is even and hence $\frac{k}{2}-u-2=2t-u-1$. Note that since $u$ is odd, $10u+6$ is divisible by $4$. Since $4\not| k$, in this subcase we see that $k \geq 10u+8$. 
 
 Since  $u+1$ and $3u+1$ both are even, we  consider the solution $(1, \ldots, 1,2t-u-2,\frac{2t+u+1}{2},\frac{2t+u+1}{2})$ and using the color of $1$, we get  $\chi(2t-u-2)=1$.  Since $\frac{k}{2}-u-2 = 2t-u-1 \geq \frac{2t+3u+1}{2} \Leftrightarrow k \geq 10u+8$, by taking the solution $ (1, \ldots, 1, 2t-3u-2,\frac{2t+3u+1}{2}, \frac{2t+3u+1}{2})$, 
  we conclude that $\chi(2t-3u-2)=1$.  Therefore  the solution $(\underbrace{1,\ldots,1}_{(k-4=4t-2)-\mbox{times}},2t-3u-2,2t-u-2,2t-u-1,2t-u-1)$ is a $2$-zero-sum solution to the equation $\mathcal{E}(k)$ irrespective of the color of $2t-u-1$.

\bigskip

\noindent{\sf Case III:} ($u$ is even)

\bigskip

Let $\chi : [1,\frac{k}{2}-u-2] \rightarrow \{0,1\}$ be a $2$-coloring such that  $\chi$ doesn't admit a $2$-zero sum solution to $\mathcal{E}(k)$. We may assume that $\chi(1) = 0$.

\bigskip

\noindent{\sf Subcase 1.} ($ 4 | k$) 

\bigskip

Since $u$ is even, $4$ divides $10u$ and hence $k \geq 10u +8$. Now, using $\chi(1)=0$ and by considering the solution $(1, \ldots, 1,\frac{k}{4}+\frac{u}{2}, \frac{k}{4}+\frac{u}{2}, \frac{k}{2}-u-3)$,  we conclude that $\chi(\frac{k}{2}-u-3)=1$. Also, since $ \frac{k}{2}-u-2 \geq \frac{k}{4}+\frac{3u}{2} \Leftrightarrow k \geq 10u+8$,  from the solution $(1, \ldots, 1, \frac{k}{4}+\frac{3u}{2}, \frac{k}{4}+\frac{3u}{2}, \frac{k}{2}-3u-3)$, we get $\chi(\frac{k}{2}-3u-3)=1$. Note that whatever color of $\frac{k}{2}-u-2$ may be, the solution $(\underbrace{1,\ldots,1}_{\text{(k-4)-times}},\frac{k}{2}-3u-3,\frac{k}{2}-u-3,\frac{k}{2}-u-2,\frac{k}{2}-u-2)$ satisfies $2$-zero sum condition to the equation $ \mathcal{E}(k)$ and finishes this subcase.

\bigskip

\noindent{\sf Subcase 2.} ($4 \not|\ k$)

\bigskip

 In this case, since $k$ is even, we can write $k =4t+2$  for some integer $t$ and $\frac{k}{2}-u-2 = 2t-u-1$. Now, consider the solution $ (1, \ldots, 1,2t-u-1,\frac{2t+u}{2},\frac{2t+u}{2})$ and using the color of $1$, we conclude $\chi(2t-u-1)=1$. Since $ \frac{k}{2}-u-2= 2t-u-1 \geq \frac{2t+3u}{2} \Leftrightarrow k \geq 10u+6$, the solution $ (1, \ldots, 1, 2t-3u-1,\frac{2t+3u}{2},\frac{2t+3u}{2})$ implies that $\chi(2t-3u-1)=1$. Then, the solution $(\underbrace{1,\ldots,1}_{(k-4=4t-2)-\mbox{times}}, 2t-3u-1, 2t-u-1, 2t-u-1, 2t-u-1)$ is a $2$-zero-sum solution to the equation $\mathcal{E}(k)$. Thus this proves the subcase and the upper bound.

\bigskip

Now we prove the lower bound for all positive integer $u$ with $k\geq 16$.

\bigskip

\noindent{\sf Case I}: ($4 | k$) 

\bigskip

\noindent{\sf Subcase 1.}  ($u$ is odd)

\bigskip

Let $\chi : [1, \frac{k}{2}- u-3] \rightarrow \lbrace 0,1\rbrace$ be a $2$-coloring with $\chi(1) = \chi(2)=\cdots=\chi( \frac{k}{4}- \frac{u+3}{2}) = 0$ and $\chi(\frac{k}{4}- \frac{u+3}{2}+1) =\cdots=  \chi(\frac{k}{2}- u-3) = 1$.

\smallskip

If there exists a solution and the color of all the $x_i$'s are $0$, then $x_1+\cdots+x_{k-3} \geq k-3$ and $x_{k-2}+x_{k-1}+x_k \leq 3(\frac{k}{4}- \frac{u+3}{2})$.  Since $k-3 >  3(\frac{k}{4}- \frac{u+3}{2})$, this is not possible.  Therefore, any  $2$-zero sum solution, at least two of $x_{k-2}, x_{k-1}$ and $x_k$ must have color $1$ and the other has color $0$ or all three have color $1$. In the first case,  we get 
$$x_{k-2}+x_{k-1}+x_k \leq \frac{k}{4}- \frac{u+3}{2}+ 2 \left(\frac{k}{2}- u-3\right) = \frac{5k-10u-30}{4}
$$  
whereas $x_1+\cdots+x_{k-3} \geq k-3$. Note that $\frac{5k-10u-30}{4} < k-3 \Leftrightarrow k < 10u+18$.  Since $4|k$, we get, 
$x_1+\cdots+x_{k-3} > x_{k-2}+x_{k-1}+x_k$ if and only if  $k \leq 10u+14 $. Therefore, in this case, we conclude that there is no $2$-zero-sum solution.

\smallskip

For the second case, we assume $\chi(x_i) = 1$ for all $i = k-2, k-1, k$ and  $\chi(x_j) = 1$ for some $j\in \{1, 2, \ldots, k-3\}$.  Hence, we get, 
$x_{k-2}+x_{k-1}+x_k \leq  3(\frac{k}{2}- u-3)$ whereas   $x_1+\cdots+x_{k-3} \geq  k-4+ (\frac{k}{4}- \frac{u+3}{2}+1) = \frac{5k-2u-18}{4}$. Note that  $3(\frac{k}{2}- u-3) < \frac{5k-2u-18}{4} \Leftrightarrow k < 10u+18 \Leftrightarrow k \leq 10u+14$.  This condition implies that the chosen color has no $2$-zero sum solution.

\bigskip

\noindent {\sf Subcase 2.} ($u$ is even)

\bigskip

In this case, consider the 2-coloring $\chi : [1, \frac{k}{2}- u-3  ] \rightarrow \lbrace 0,1\rbrace$ with $\chi(1) = \cdots = \chi(\frac{k}{4}- \frac{u+2}{2}) = 0 $ and $\chi(\frac{k}{4}- \frac{u+2}{2}+1) = \cdots =  \chi(\frac{k}{2}- u-3) = 1$.

\smallskip

If there exists a solution and the color of all the $x_i$'s are $0$, then $x_1+\cdots+x_{k-3} \geq k-3$ and $x_{k-2}+x_{k-1}+x_k \leq 3(\frac{k}{4}- \frac{u+2}{2})$.  Since $k-3 >  3(\frac{k}{4}- \frac{u+2}{2})$, this type of solution is not possible.  Therefore, any  $2$-zero sum solution, at least two of $x_{k-2}, x_{k-1}$ and $x_k$ must have color $1$ and the other has color $0$ or all three have color $1$. Hence, in the first case, we get 
$$x_{k-2}+x_{k-1}+x_k \leq \frac{k}{4}- \frac{u+2}{2}+ 2 \left(\frac{k}{2}- u-3\right) = \frac{5k-10u-28}{4}
$$  
whereas $x_1+\cdots+x_{k-3} \geq k-3$. Note that $\frac{5k-10u-28}{4} < k-3 \Leftrightarrow k < 10u+16$.  Since $4|k$, we get, 
$x_1+\cdots+x_{k-3} > x_{k-2}+x_{k-1}+x_k$ if and only if  $k \leq 10u+14 $. Therefore, in this case, we conclude that there is no $2$-zero-sum solution.

\smallskip

For the second case, we assume $\chi(x_i) = 1$ for all $i = k-2, k-1, k$ and  $\chi(x_j) = 1$ for some $j\in \{1, 2, \ldots, k-3\}$.  Hence, we get, 
$x_{k-2}+x_{k-1}+x_k \leq  3(\frac{k}{2}- u-3)$ whereas   $x_1+\cdots+x_{k-3} \geq  k-4+ (\frac{k}{4}- \frac{u+2}{2}+1) = \frac{5k-2u-16}{4}$. Note that  $3(\frac{k}{2}- u-3) < \frac{5k-2u-16}{4} \Leftrightarrow k < 10u+20 \Leftrightarrow k \leq 10u+18$. Thus, for $k\leq 10u +14$, we conclude that  the chosen color has no $2$-zero sum solution. This finishes the proof of the lower bound.

\bigskip

\noindent{\sf Case II:} ($4 \not|\ k$)

\bigskip

 In this case,  we can write $k= 4t +2 $ for some integer $t$. 

\bigskip

\noindent{\sf Subcase 1.}  ($u$ is even)

\bigskip

 For proving the lower bound, consider the $2$-coloring $\chi : [1, \frac{4t+2}{2}- u-3] \rightarrow \lbrace 0,1\rbrace$ with $\chi(1) = \cdots= \chi(t-\frac{u+2}{2}) = 0$ and $\chi(t-\frac{u+2}{2}+1) = \cdots=\chi(2t-u-2) = 1$.

\smallskip

If there exists a solution and the color of all the $x_i$'s are $0$, then $x_1+\cdots+x_{k-3} \geq k-3 =4t-1$ and $x_{k-2}+x_{k-1}+x_k \leq 3(t- \frac{u+2}{2})$.  Since $4t-1 >  3(t- \frac{u+2}{2})$, this type of solution is not possible.  Therefore, in any  $2$-zero sum solution, at least two of $x_{k-2}, x_{k-1}$ and $x_k$ must have color $1$ and the other has color $0$ or all three have color $1$. Hence, in the first case, we get 
$$x_{k-2}+x_{k-1}+x_k \leq t- \frac{u+2}{2}+ 2 \left(2t- u-2\right) = \frac{10t-5u-10}{2}
$$  
whereas $x_1+\cdots+x_{k-3} \geq k-3 = 4t-1$. Note that $\frac{10t-5u-10}{2} < 4t-1 \Leftrightarrow k < 10u+18$.  Since $4\not|k$, we get, 
$x_1+\cdots+x_{k-3} > x_{k-2}+x_{k-1}+x_k$ if and only if  $k \leq 10u+14 $. Therefore, in this case, we conclude that there is no $2$-zero-sum solution.

\smallskip

For the second case, we assume $\chi(x_i) = 1$ for all $i = k-2, k-1, k$ and  $\chi(x_j) = 1$ for some $j\in \{1, 2, \ldots, k-3\}$. Hence, we get, 
$x_{k-2}+x_{k-1}+x_k \leq  3(2t- u-2)$ whereas   $x_1+\cdots+x_{k-3} \geq  k-4+ (t- \frac{u+2}{2}+1) = \frac{10t-u-4}{2}$. Note that  $3(2t- u-2) < \frac{10t-u-4}{2} \Leftrightarrow k < 10u+18 \Leftrightarrow k \leq 10u+16$. Hence, for $k\leq 10u +14$, we conclude that  the chosen color has no $2$-zero sum solution. 

\bigskip

\noindent{\sf Subcase 2.}  ($u$ is odd)

\bigskip

 For proving the lower bound, consider the $2$-coloring $\chi : [1, \frac{4t+2}{2}- u-3] \rightarrow \lbrace 0,1\rbrace$ with $\chi(1) = \cdots = \chi(t-\frac{u+3}{2})= 0$ and $\chi(t-\frac{u+3}{2}+1) = \cdots = \chi(2t-u-2) = 1$.

\smallskip

If there exists a solution and the color of all the $x_i$'s are $0$, then $x_1+\cdots+x_{k-3} \geq k-3 =4t-1$ and $x_{k-2}+x_{k-1}+x_k \leq 3(t- \frac{u+3}{2})$.  Since $4t-1 >  3(t- \frac{u+3}{2})$, this type of solution is not possible.  Therefore, in any  $2$-zero sum solution, at least two of $x_{k-2}, x_{k-1}$ and $x_k$ must have color $1$ and the other has color $0$ or all three have color $1$. Hence, in the first case, we get 
$$x_{k-2}+x_{k-1}+x_k \leq t- \frac{u+3}{2}+ 2 \left(2t- u-2\right) = \frac{10t-5u-11}{2}
$$  
whereas $x_1+\cdots+x_{k-3} \geq k-3 = 4t-1$. Note that $\frac{10t-5u-11}{2} < 4t-1 \Leftrightarrow k \leq 10u+18$.  Since $4\not|k$, we get, 
$x_1+\cdots+x_{k-3} > x_{k-2}+x_{k-1}+x_k$ if and only if  $k \leq 10u+16 $. Therefore, in this case, we conclude that there is no $2$-zero-sum solution.

\smallskip

For the second case, we assume $\chi(x_i) = 1$ for all $i = k-2, k-1, k$ and  $\chi(x_j) = 1$ for some $j\in \{1, 2, \ldots, k-3\}$.  Hence, we get, 
$x_{k-2}+x_{k-1}+x_k \leq  3(2t- u-2)$ whereas   $x_1+\cdots+x_{k-3} \geq  k-4+ (t- \frac{u+3}{2}+1) = \frac{10t-u-5}{2}$. Note that  $3(2t- u-2) < \frac{10t-u-5}{2} \Leftrightarrow k < 10u+16.$ Thus, for $k\leq 10u +14$, we conclude that  the chosen color has no $2$-zero sum solution. This proves the lower bound and the theorem. \hfill\(\Box\)

\section{proof of the theorem \ref{more}}

Here we consider the equation
 \begin{equation}\label{eqmore}
 \mathcal{E}_{v}^{(1,0)}:x_1 + \cdots +x_{k-vr} = x_{k-vr+1} + \cdots +x_k.
 \end{equation}

Let us denote $s:= \frac{k}{r} - \left\lfloor \frac{(v-1)k}{vr} \right\rfloor -1$  for simplicity. 

\medskip
 
\noindent{\sf Case I:}  ($k=2vr$) 

\medskip

In this case, the number of variables in both sides of the \eqref{eqmore} are equal and hence $(1, \ldots, 1)$ is an $r$-zero sum solution.
Thus we get, $S_{\mathfrak{z},r}^{(1,0)}(k;r;v) = 1$, as desired. 

\medskip

\noindent{\sf Case II:} ($k > 2vr$)

\medskip

Since $r|k$, by division algorithm we write $k-2vr = vrt + ir$ for some non-negative integers $t$ and  $i$ with $i\in [1, v]$. Therefore, we get 
\begin{align*}
k  = vrt +(2v+i)r &\Leftrightarrow(v-1) k  = v(v-1)rt +(2v+i)(v-1)r \\
& \Leftrightarrow \frac{(v-1)k}{vr}  = (v-1)t +\frac{(v-1)(2v+i)}{v}\\
& \Leftrightarrow \left\lfloor \frac{(v-1)k}{vr} \right\rfloor  = (v-1)t +(2v+i-3).
\end{align*}
Hence, we have $ s = \frac{k}{r} - \left\lfloor \frac{(v-1)k}{vr} \right\rfloor -1 = \frac{k}{r}-(v-1)t -(2v+i-3)-1.$
\bigskip

For proving the lower bound we show that  in the interval $[1, s-1]$, equation \eqref{eqmore} does not have any solution. First, we observe that
\begin{align*}
&x_{k-vr+1} + \cdots +x_k  \leq  vr(s-1)  =  vr \left(\frac{k}{r} - \left\lfloor \frac{(v-1)k}{vr} \right\rfloor -1 -1\right) \\
                      & =  vr \left(\frac{k}{r} - (v-1)l -(2v+i-3)-1 -1\right)=  vk- v(v-1)rl-v(2v+i-3)r-vr-vr\\
                      & =  vk -(v-1)k +(2v+i)(v-1)r-v(2v+i-3)r -2vr \\
                      & =  k-(2v+i)r+3vr-2vr =k -vr -ir
\end{align*} but  $x_1 + \cdots +x_{k-vr} \geq k -vr$. Thus, we get $S_{\mathfrak{z},r}^{(1,0)}(k;r;v) \geq  \frac{k}{r} - \left\lfloor \frac{(v-1)k}{vr} \right\rfloor -1$.

\bigskip

Now, we proceed to prove the upper bound and consider an arbitrary $m$-coloring $ \chi :[1, s] \rightarrow \{0,1, \ldots , m-1\}$. Next, we consider 
\smallskip
\begin{equation}\label{solution}
 (x_1,x_2, \ldots, x_{k-vr}, x_{k-vr+1}, \ldots, x_k) = (\underbrace{1,\ldots,1}_{(k-vr)-\text{times}} ,\underbrace{s,\ldots,s}_{ir-\text{times}}, \underbrace{s-1, \ldots, s-1}_{(vr-ir)-{\text times}})
\end{equation}
 and show that it is a solution to the equation \eqref{eqmore}. Note that the value of $x_1 + \cdots +x_{k-vr}$ is $ k-vr$ and the value of $x_{k-vr+1}+ \ldots +x_k$ is $irs +(v-i)r (s-1)$ which is exactly equal to $k-vr$ because,
\begin{align*}
 &irs +(v-i)r (s-1)=ir\left(\frac{k}{r} - \left\lfloor \frac{(v-1)k}{vr} \right\rfloor -1\right) +(v-i)r\left(\frac{k}{r} - \left\lfloor \frac{(v-1)k}{vr} \right\rfloor -1-1\right)\\
                        & = ir\left(\frac{k}{r} - (v-1)l -(2v+i-3)-1)+(v-i\right)r\left(\frac{k}{r} - (v-1)l -(2v+i-3)-1-1\right)\\
                        & = vk-v(v-1)rl-(2v+i-3)vr-2vr+ir\\
                        & = vk-(v-1)k +(2v+i)(v-1)r-(2v+i-3)vr-2vr+ir = k-vr.                                \end{align*}

\smallskip

Thus, it remains to show that the tuple defined in \eqref{solution} satisfies the $r$-zero sum condition. Since $r \mid k$, we see that
\begin{align*}
\sum_{i=1} ^k \chi (x_i) & = (k-vr) \chi (1) + ir ( \chi (s)) + (vr-ir) \chi (s-1) \equiv 0 (\textrm{mod}\ r).
\end{align*} Therefore, we get $S_{\mathfrak{z},r}^{(1,0)}(k;r;v) \leq  \frac{k}{r} - \left\lfloor \frac{(v-1)k}{vr} \right\rfloor -1$ and hence the theorem.
$\hfill\Box$

\vskip 20pt
\begin{acknowledgement*}  
We would like to sincerely thank Prof. S. D. Adhikari and Prof. R. Thangadurai for their insightful remarks. This work was done while the authors were visiting the Department of Mathematics, Ramakrishna Mission Vivekananda Educational and Research Institute and they wish to thank this institute for the excellent environment and for their  hospitality.
\end{acknowledgement*}

\end{document}